\documentclass[12pt,reqno]{amsart}

\numberwithin{equation}{section}

\usepackage[final, 
color]{showkeys}
\definecolor{refkey}{gray}{.85}
\definecolor{labelkey}{gray}{.85}

\usepackage{AKstyle}

\usepackage{enumerate}

\usepackage{caption}
\usepackage[labelfont=rm]{subcaption}

\usepackage[pagebackref=true, colorlinks]{hyperref}

\hypersetup{pdffitwindow=true,linkcolor=blue,citecolor=blue,urlcolor=blue,menucolor=blue}

\usepackage{comment}

\newcommand{\bbF}{\mathbb F}
\renewcommand{\F}{\cQ}

\begin{document}

\author{Jiuzu Hong}
%\thanks{JH is partially supported by...}
\email{jiuzu.hong@yale.edu}
\address{Department of Mathematics, Yale University, New Haven, CT}

\author{Alex Kontorovich}
\thanks{A.K. is partially supported by
an NSF CAREER grant DMS-1254788, an Alfred P. Sloan Research Fellowship, a Yale Junior Faculty Fellowship, and support at IAS from The Fund for Math and The Simonyi Fund.}
\email{alex.kontorovich@yale.edu}
\address{Department of Mathematics, Yale University, New Haven, CT, and School of Mathematics, IAS, Princeton, NJ}

\title[Anisotropic and Thin Pythagorean Orbits]{Almost Prime Coordinates for Anisotropic and Thin Pythagorean Orbits}
%A Remark on the Affine Sieve}

\begin{abstract}
We make an  observation which doubles the exponent of distribution 
in certain Affine Sieve problems, such as those considered by Liu-Sarnak, Kontorovich, and Kontorovich-Oh. 
As a consequence, we decrease the known bounds on the saturation numbers in these problems.
%produce $R$-almost primes in these settings, with much smaller values of $R$.
\end{abstract}
\date{\today}
\subjclass[2010]{}
\maketitle
\tableofcontents

\section{Introduction}\label{sec:intro}

The purpose of this paper is to make a simple  observation about the execution of the Affine Sieve, which has the effect of doubling the exponent of distribution in many natural  sieve problems. For concreteness, we will illustrate the method on problems studied by Liu-Sarnak \cite{LiuSarnak2010}, Kontorovich \cite{MyThesis, Kontorovich2009}, and  Kontorovich-Oh \cite{KontorovichOh2012}. We first briefly recall the general %procedure
setup, then specialize to these particular problems,  explain what is known, and finally describe our results. % innovation.

\subsection{The General Affine Sieve}\

Roughly speaking, the Affine Sieve inputs a pair $(\cO,f)$ consisting of $(i)$ an integer orbit $\cO\subset\Z^{n}$ by a linear group and $(ii)$ a polynomial function $f$ which is integral on $\cO$, and outputs a number $R\le\infty$ so that there are ``many'' points $\bx\in\cO$ with $f(\bx)$ having at most $R$ prime divisors. Let us make this precise. 

Let $\G<\GL_{n}(\Z)$ be a finitely-generated group of invertible $n\times n$ integer matrices, let $\bbG:=\Zcl(\G)$ be its Zariski closure, and 
denote the real points of
% its %(affine) 
%Zariski closure. 
$\bbG$
by $G:=\bbG(\R)$.
When the Haar measure of $\G\bk 
G$ is infinite, we refer to $\G$ as {\it thin}. For a fixed primitive vector $\by\in\Z^{n}$, we consider the orbit 
\be\label{eq:cO}
\cO\ :=\ \by\cdot \G \ \subset \ \Z^{n};
\ee
we refer to $\cO$ as thin when $\G$ is. 
Given a polynomial $f
%:\R^{n}\to\R
$ in $n$ variables which is integral on $\cO$,
 we say that the pair $(\cO,f)$ is {\it strongly primitive}%
\footnote{%
This corrects the
 %definition given 
 typo
 in \cite[Definition 1.3]{KontorovichOh2012}.%
} 
if,
for all integers $q\ge1$, there is an $\bx\in\cO$ so that $f(\bx)\in(\Z/q\Z)^{\times}$. We assume henceforth that this is the case.

For an integer $R\ge1$, let 
$$
\cP_{R}\ \subset\ \Z
$$ 
denote the set of $R$-almost primes, that is, numbers with at most $R$ prime divisors. We allow $R=\infty$, in which case $\cP_{R}=\Z$.
%
%The Affine Sieve is concerned with producing ``many'' points $\bx\in\cO$ 
%
%Given a pair $(\cO,f)$ as above and any $R\le\infty$, l
Let
$$
\cO(f,R)\ :=\ 
\{\bx\in\cO\ : \ f(\bx)\in\cP_{R}\}
.
$$
The goal of the Affine Sieve is
% the 
to
study % of
 the so-called {\it saturation number}, 
%defined as
given by
$$
R_{0}(\cO,f)
\
:=\ \min\{R\le\infty\ : \ \Zcl(\cO(f,R))=\Zcl(\cO)\}
.
%,
$$
That is, $R_{0}$ is  the minimal $R$ for which 
$\cO(f,R)$
is Zariski dense in the Zariski closure of $\cO$. (Here $\Zcl(\cdot)$ refers to the Zariski closure in affine  space $\bA_{\Q}^{n}$.)

The program 
initiated  by Bourgain-Gamburd-Sarnak \cite{BourgainGamburdSarnak2006, BourgainGamburdSarnak2010} and completed  by Salehi Golsefidy-Sarnak \cite{SalehiSarnak2011}
shows in essentially  the greatest generality possible that pairs $(\cO,f)$ are {\it factor finite}, meaning that we have the strict inequality $R_{0}(\cO,f)<\infty$. 
Beyond factor-finiteness, one would like to actually determine the saturation number of any given pair $(\cO,f)$.
As stated, this problem is completely hopeless, as it includes 
all classical sieve problems (see the discussion in, e.g., \cite{BourgainGamburdSarnak2010}).
Nevertheless, there is an ongoing program of determining, or at least giving strong estimates for, the saturation number in certain specific cases, where more structure can be exploited.
%Most sieve-type  problems (including the classical ones) %boil down 
%amount
%to % bounding 
%determining
%$R_{0}$ for various choices of $\cO$ and $f$, 
We give a few natural examples below.

\subsection{Thin  Pythagorean Orbits}\label{sec:thinPyth}\

Let $G=\SO_{\F}(\R)<\SL_{3}(\R)$ be the real special orthogonal group preserving the ``Pythagorean'' quadratic form
\be\label{eq:PythForm}
\F(\bx) %=\F_{P}(\bx)
\ :=\ x^{2}+y^{2}-z^{2}
,
\ee 
where $\bx=(x,y,z)$.
We define a Pythagorean triple to be a primitive integer vector on the cone $\F=0$.

 Let $\G<G(\Z)$ be a finitely generated subgroup of the integer matrices in $G$, and assume $\G$ is non-elementary, 
 %that is, 
 or equivalently, that
 its Zariski closure is $\SO_{\F}$. For a fixed Pythagorean triple, e.g., $\by=(3,4,5)$, let $\cO$ be its  corresponding $\G$-orbit, as in \eqref{eq:cO}. We allow $\G$, and hence $\cO$, to be thin, in which case we refer to $\cO$ as a thin Pythagorean orbit.

A measure of this thinness is the critical exponent 
$$
\gd=\gd_{\G} \ \in \ [0,1]
$$ of $\G$; this is the abscissa of convergence of the Poincar\'e series of $\G$, or equivalently, the Hausdorff dimension of the limit set of $\G$.
Since $\G$ is non-elementary, $\gd$ is strictly positive; moreover $\G$ is thin if and only if $\gd<1$.
The role played by this geometric invariant is illustrated by the easy 
fact that
\be\label{eq:cOsize}
\#\{\bx\in\cO\ : \ \|\bx\|<T\} \ = \ T^{\gd+o(1)},
\ee
as $T\to\infty$, where $\|\cdot\|$ is the standard Euclidean norm. 

For various choices of the polynomial $f$, one can consider the problem of estimating the saturation number $R_{0}(\cO,f)$. Three natural choices for $f$ considered in \cite{MyThesis, Kontorovich2009, KontorovichOh2012} are
\be\label{eq:fsPyth}
\threecase
{f_{\sH}(\bx)=z,}{the ``hypotenuse'',}
{f_{\sA}(\bx)=\frac1{12}xy,}{the ``area'',}
{f_{\sC}(\bx)=\frac1{60}xyz,}{the product of coordinates.}
\ee
Recall that we assume, as throughout, that the pair $(\cO,f)$ is strongly primitive; the fractions in \eqref{eq:fsPyth} are to remove extraneous prime factors (e.g. it is elementary that the product of coordinates $xyz$  in a Pythagorean triple is always divisible by $60$).

We will refer to the pairs $(\cO,f)$ above with $f\in\{f_{\sH},f_{\sA},f_{\sC}\}$ as Examples $A$, $B$, and $C$, respectively.

\begin{thm}[\cite{MyThesis, Kontorovich2009, KontorovichOh2012}]\label{thm:KO}
Assume the critical exponent $\gd$ of $\G$ is sufficiently close to $1$. Then % for 
%$f\in\{f_{\sH},f_{\sA},f_{\sC}\}$, 
we have 
\be\label{eq:Rpyth}
R_{0}(\cO,f) \ \le \
\threecase
{13,}{for Example $A$,}%{if $f=f_{\sH}$,}
{40,}{for Example $B$,}%{if $f=f_{\sA}$,}
{58,}{for Example $C$.}%{if $f=f_{\sC}$.}
\ee
\end{thm}
\begin{rmk}\label{rmk:KO}
 We have taken this opportunity to correct
the values of $R$ 
 in the statement of \cite[Theorem 1.5]{KontorovichOh2012}, which 
 were improperly computed; see Remark \ref{rmk:corrKO}.
 %the footnote on p. \pageref{foot}.
 % \S\ref{sec:?} for the corrected procedure, used to produce the numbers in \eqref{eq:Rpyth}.
 % (the correct value of $\mu$ in \cite[(5.9)]{KontorovichOh2012} is given in \cite{eq:muIs} below).
\end{rmk}
 
% \newpage

% in the above to correct the numerical values of $R$ claimed in \cite{KontorovichOh2012}, which were inaccurately computed. 

The upper bounds on $R_{0}$ given in \eqref{eq:Rpyth} are based on Gamburd's spectral gap \cite{Gamburd2002} (see \S\ref{sec:spec}); the lower bounds, and expected true values of $R_{0}$, are the so-called ``sieve dimensions'' %$\gk$ 
(see 
%\S\ref{sec:sieve}
Remark \ref{rmk:dim}), given by
\be\label{eq:gkPyth}
\gk=\gk(\cO,f)\ := \
\threecase
{1,}{for Example $A$,}%{if $f=f_{\sH}$,}
{4,}{for Example $B$,}%{if $f=f_{\sA}$,}
{5,}{for Example $C$.}%{if $f=f_{\sC}$.}
\ee

In % the ``hypotenuse'' case of $f=f_{\sH}$
Example $A$, the upper bound on the saturation number has been reduced significantly  in \cite{BourgainKontorovich2013} to 
\be\label{eq:BK}
R_{0}(\cO,f_{\sH})\ \le \ 4
\ee
 by
 % completely 
 quite
 different methods from those discussed here, so we will not focus on this case.
  %here. 
  For the other two choices of $f$, an easy consequence of our method is the following improvement.

\begin{thm}\label{thm:Pyth}
Theorem \ref{thm:KO} holds with \eqref{eq:Rpyth} replaced by
\be\label{eq:RpythNew}
R_{0}(\cO,f) \ \le \
\twocase{}
{25,}{for Example $B$,}%{if $f=f_{\sA}$,}
{37,}{for Example $C$.}%{if $f=f_{\sC}$.}
\ee
\end{thm}

\begin{rmk}\label{rmk:ExA}
For Example $A$, our method gives $R_{0}(\cO,f_{\sH})\le 7$; see \S\ref{sec:ExA}. This is an improvement over \eqref{eq:Rpyth}, but does not compete with \eqref{eq:BK}.
\end{rmk}

\begin{rmk}
In all the statements above  (and below), the Zariski density is an easy consequence of a lower bound on the cardinality of $\cO(f,R)$ restricted roughly to an archimedean ball, see
Remark \ref{rmk:Zariski}.
% \eqref{eq:anTLower}.
\end{rmk}

%Other consequences are discussed in \S\ref{sec:?}. 

\subsection{Anisotropic Orbits}\label{sec:ani}\

In \cite{LiuSarnak2010}, Liu-Sarnak consider a related problem. Instead of the isotropic Pythagorean form $\F$ in \eqref{eq:PythForm}, they let $\F$ be an {\it anistropic} (over $\Q$) indefinite integral ternary quadratic form, e.g. $\F(\bx)=x^{2}+y^{2}-3z^{2}$. This means that there are no rational points on the cone $\F=0$, and
so   one instead  considers the affine quadric
\be\label{eq:levelSet}
 V=V_{\F,t}\ := \ \{\bx: \F(\bx)=t\},
\ee
 for a fixed non-zero integer $t$, chosen so that $V(\Z)$ is non-empty. 
 For simplicity, assume that  $t\cdot\gD(\F)$ is square-free, where $\gD(\F)$ is the discriminant of $\F$. 
The study of the vectors in $V(\Z)$ reduces (see \cite[\S2]{LiuSarnak2010}) to that of orbits $\cO:=\by\cdot \G$, where
$\by\in V(\Z)$, and
 $\G=\SO_{\F}(\Z)$ is the integer matrix group preserving $\F$. (Such an orbit is not thin, as $\G$ is a lattice in $G=\SO_{\F}(\R)$, with critical exponent $\gd=1$.)  Let $f(\bx)=xyz$ be the product of coordinates, and recall our assumption that the pair $(\cO,f)$ is strongly primitive (for example, this is guaranteed if $\by=(1,1,1)$).

We refer to this pair $(\cO,f)$ as Example $D$.
 
 \begin{thm}[\cite{LiuSarnak2010}]\label{thm:LS}
%With notation as above,  w
We have the following bound on the saturation number in Example $D$:
\be\label{eq:RLS}
R_{0}(\cO,f)\ \le \ 26.
\ee
Assuming the Selberg Eigenvalue Conjecture (see Theorem \ref{thm:spec}), we have
\be\label{eq:RLSS}
R_{0}(\cO,f)\ \le \ 22.
\ee
 \end{thm}
 
% We will denote by Example $E$ the same pair $(\cO,f)$ but under the further assumption of the Selberg Eigenvalue Conjecture.

\begin{rmk}\label{rmk:IsoAni0}
Note that for the product of coordinates here the sieve dimension is $\gk=3$ (see \S\ref{sec:sieve}), rather than $\gk=5$ 
in \eqref{eq:gkPyth} for Example $C$, that is, for $f=f_{\sC}$. This is because in the isotropic case, there are 
 non-constant polynomial parametrizations of the integer points of the corresponding orbits which can be (and, in the Pythagorean case, are) reducible; see Remark \ref{rmk:IsoAni}. 
\end{rmk}

As a consequence of our method, we have
\begin{thm}\label{thm:Aniso}
Theorem \ref{thm:LS} holds unconditionally %(that is, for Example $D$) 
with \eqref{eq:RLS} replaced by 
$$
R_{0}(\cO,f)\ \le \ 16.
$$ 
Assuming the Selberg Eigenvalue Conjecture% (that is, Example $E$)
, \eqref{eq:RLSS} may be replaced by 
$$
R_{0}(\cO,f)\ \le \ 14.
$$ 
\end{thm}

\subsection{New Observation}\label{sec:new}\

Our key new observation is that, for all the problems above (indeed for nearly all natural  Affine Sieve problems in the literature), the polynomial $f$ is %actually 
{\it homogeneous}. Roughly speaking, this allows us, in the modular/archimedean decomposition of the Affine Sieve, to projectivize, taking a larger stabilizer group (see \S\ref{sec:obs}). %This leads to
As a result, we have
 no modular loss in the error terms, whereas in the previous approaches, a power of the level was lost; see Remark \ref{rmk:new}.
%. We explain this statement more precisely in \S\ref{sec:affStd}. 
The upshot is an improvement by a factor of two in the level of distribution (see \S\ref{sec:level1}) in the above problems, which translates to the above-claimed improved bounds on saturation numbers.
%
%\begin{rmk}
In fact our main observation is a general principle, applying to many other settings, e.g., the pairs $(\cO,f)$ considered in \cite{NevoSarnak2010} with $f$ homogeneous; we will not %further pursue 
bother with
other applications here. 
%\end{rmk}

\subsection{Outline}\
In \S\ref{sec:2} we collect some relevant background. In particular, we recall facts on spectral gaps, counting, levels of distribution, and the Diamond-Halberstam-Richert sieve. We also  sketch proofs of Theorems \ref{thm:KO} and \ref{thm:LS}, since our proofs of Theorems \ref{thm:Pyth} and \ref{thm:Aniso} are nearly identical. In \S\ref{sec:pfs} we explain our new observation, and use it to prove Theorems  \ref{thm:Pyth} and \ref{thm:Aniso}. %We conclude in \S\ref{sec:concl} with some comments on future directions.

\subsection{Notation}\

We use the standard notation $f=O(g)$ and  $f\ll g$ synonymously  to mean $f(x)\le C g(x)$ for an implied constant $C>0$ and all $x$ sufficiently large. Unless otherwise specified, $C$ may depend only on the pair $(\cO,f)$, which is treated as fixed. The little-oh notation $f=o(g)$ means $f/g\to0$.

\subsection*{Acknowledgements}\

The authors thank Shamgar Gurevich, Nick Katz, and Peter Sarnak for enlightening discussions.

\newpage

\section{Background}\label{sec:2}

\subsection{Spectral Gap}\label{sec:spec}\

We take the following as our definition of a spectral gap for the cases of interest to us here. 
For $\cQ$ a ternary indefinite integral quadratic form (either isotropic or anisotropic over $\Q$), let
 $$
 G=\SO_{\cQ}(\R)\cong \SO_{2,1}(\R)
 $$
 be its stabilizer group,
  and let $\G<G(\Z)$ be a finitely generated (and hence geometrically finite)  integer 
% (and hence discrete) 
subgroup with critical exponent 
$$
\gd>1/2.
$$ 
The decomposition of the right regular representation of $G$ on 
$L^{2}(\G\bk G)$ is of the form \cite{GelfandGraevPS1966, LaxPhillips1982}
$$
L^{2}(\G\bk G)
=
V_{0}\oplus V_{%\vf_
{1}}\oplus\cdots\oplus V_{%\vf_
{J}}\oplus V_{temp}.
$$
Here
 $V_{temp}$ is a (reducible) subspace consisting  of the tempered spectrum;
 the $V_{%\vf_
 {j}}$, $j=1,\dots,J$ are %unitary, 
 %irreducible, and 
 isomorphic as $G$-representations to complementary series representations with corresponding parameters 
 $$
 1/2<s_{J}\le\cdots\le s_{1}<\gd\le 1
 $$ 
 (in our normalization, the principal series representations lie on the critical line $\Re(s)=1/2$);
and
 $V_{0}$ is either the trivial representation  if $\G$ is a lattice, or a complementary series representation of parameter $s_{0}=\gd$ if $\gd<1$ \cite{Patterson1976, Sullivan1984}. We say a number $s\in(1/2,1)$ {\it appears} in $L^{2}(\G\bk G)$ if it arises as one of the $s_{j}$ above.

For a square-free integer $q\ge1$, define the  level $q$ principal congruence subgroup of $\G$ as
$$
\G(q)\ :=\ \{\g\in\G\ : \ \g\equiv I(\mod q)\}.
$$
We have a similar decomposition for $L^{2}(\G(q)\bk G)$, and the inclusion $\G(q)<\G$ induces the reverse inclusion on spectrum; that is, any %complementary series representation of 
parameter $s$ which appears in $L^{2}(\G\bk G)$ also appears in $L^{2}(\G(q)\bk G)$. We say $s\in(1/2,1)$ is  {\it new} spectrum at level $q$, if the parameter $s$ appears in $L^{2}(\G(q)\bk G)$ but does not arise in this way as a lift from $L^{2}(\G\bk G)$; let $\Spec^{new}(q)$ denote the new spectra at level $q$.

We say that $\G$ has a {\it uniform spectral gap} 
$$
\foh\ \le \
\gt=\gt(\G)\ 
%\in\ \left[\foh,
<\
\gd %\right)
$$
 if there exists an integer 
 \be\label{eq:fBDef}
 \fB\ge1
 \ee 
 so that, for all $q$ coprime to $\fB$, 
$$
\Spec^{new}(q)\ \subset\ (1/2,\gt].
$$
In particular, the ``base'' parameter $\gd$ remains isolated as $q$ ranges through square-free numbers coprime to the ``bad'' modulus $\fB$. 

In %the 
our
archimedean (as opposed to combinatorial, for which see \cite{SalehiVarju2012}) setting, the following is the current state of affairs 
%in our setting 
on spectral gaps. % for our context.

\begin{thm}\label{thm:spec}
Assume $\gd>1/2$. Then
\begin{itemize}
\item
 $\G$ has some spectral gap $\gt\in[1/2,\gd)$ \cite{BourgainGamburd2008, BourgainGamburdSarnak2010, BourgainGamburdSarnak2011}.
\item
If moreover $\gd>5/6$, then we can take $\gt=5/6$ \cite{Gamburd2002}.
\item
If moreover $\G$ is a congruence group (and hence $\gd=1$), then we can take $\gt=1/2+7/64=39/64$ and $\fB=1$ \cite{JacquetLanglands,  KimSarnak2003}.
\item
If moreover we assume the Selberg Eigenvalue Conjecture, then we can take  $\gt=1/2$ and $\fB=1$  \cite{Selberg1965}.
\end{itemize}
\end{thm}

\begin{rmk}
In the case $\cQ$ is anisotropic over $\Q$, that is, for Example $D$% and $E$
, the quotient $\G\bk G$ is compact, and the
Jacquet-Langlands correspondence is used to apply the best-known bounds towards the Selberg Eigenvalue (or Generalized Ramanujan) Conjecture in the statement of Theorem \ref{thm:spec}.
\end{rmk}

\subsection{Effective Counting on Congruence Towers}\

With the spectral gap in place, we state the following now-standard smooth counting theorem (see, e.g., \cite{BourgainKontorovichSarnak2010} or \cite[Theorem 2.9]{BourgainKontorovich2013}% which uses 
). We define a norm on $G$ via $\|g\|^{2}=\tr g{}^{t}g$.

\begin{thm}
Assume $\G<G$ is a finitely-generated discrete group as above with critical exponent $\gd>1/2$ and spectral gap $1/2\le\gt<\gd$. Then for $T\to\infty$, there is a function $\gU_{T}:G\to\R_{\ge0}$ with the following properties. 

$(i)$ $\gU_{T}$ is a smoothed indicator of $\|g\|<T$, in the sense that
\be\label{eq:gUTsupp}
\gU_{T}(g)\ = \
\threecase
{1,}{if $\|g\|<\foh T$,}
{0,}{if $\|g\|>2 T$,}
{\in[0,1],}{otherwise,}
\ee
and\footnote{%
%Again, 
Throughout, only the exponents will be relevant to our analysis, so we will be quite crude with such statements, even when much more information is available. %one can produce an asymptotic formula here with a power savings, but this is not needed for our applications.
}
\be\label{eq:gUTsize}
\sum_{\g\in\G}
\gU_{T}(\g)\ = \ T^{\gd+o(1)}.
\ee
Moreover,

$(ii)$ for any $\g_{0}\in\G$, any square-free $q\ge1$ coprime to $\fB$ in \eqref{eq:fBDef}, and any $\Xi(q)$ satisfying $\G(q)\leqq\Xi(q)\leqq\G$, we have
\be\label{eq:gUTAP}
\sum_{\g_{1}\in\Xi(q)}\gU_{T}(\g_{1}\g_{0}) \ = \
\frac{1}{[\G:\Xi(q)]}
\sum_{\g\in\G}\gU_{T}(\g)
\ +\
O(T^{\gt+o(1)})
.
\ee
The implied constant above is independent of $q$ and $\g_{0}$. 
\end{thm}
\begin{rmk}
The interpretation of \eqref{eq:gUTAP} is that one has effective (with power savings down to the spectral gap) equidistribution of $\G$ along %arithmetic progressions 
congruence towers
mod $q$. It is important here (in fact absolutely crucial to our observation!)  to have the flexibility to choose any $\Xi(q)$ lying between $\G$ and the full level $q$ principal congruence subgroup $\G(q)$.
\end{rmk}

\subsection{Level of Distribution}\label{sec:level1}\

%, and let $\fB\ge1$ be as in Theorem \ref{thm:spec}; in particular, $\fB=1$ if $\G$ is a congruence group. 
We now define a certain  finite
sequence 
$$
\cA
\ =\ \{a_n(T)\}
$$  %be a
 of non-negative real numbers depending on a parameter
$$
T\ \to\ \infty
,
$$
which will play a key role in the analysis. This sequence is supported on values of $f(\bx)$, with $\bx\in\cO$, where the pair
%
%Let
 $(\cO,f)$
  %with $\cO=\by\cdot\G$ be 
 is one of the pairs discussed  in \S\ref{sec:thinPyth} or \S\ref{sec:ani}, that is, Examples $A$--$D$.
For ease of exposition,
we
 assume henceforth  that $\fB=1$; minor adjustments are needed in the general case. 
Using the smooth counting function from the previous subsection,  we define:
\be\label{eq:anTIs}
a_{n}(T) \ := \
\sum_{\g\in\G}
\gU_{T}(\g)
\cdot
\bo_{\{f(\by\cdot\g)=%\pm 
n\}}
.
\ee
Thus $a_{n}(T)$ counts roughly the number of representations of $n$ of the form $f(\by\cdot \g)$, for $\g$ restricted to an archimedean ball.
(In the case that $\by$ has a non-trivial stabilizer in $\G$, this will be an over-count; %nevertheless,   
statements about the Zariski closure of $\cO(f,R)$ are not sensitive to this over-counting.)

We first determine the total amount of ``mass'' contained in $\cA$, 
that is, we have from \eqref{eq:gUTsize} the approximation  
\be\label{eq:cAsize}
|\cA| \ := \ \sum_{n} a_{n}(T)  \ =\  T^{\gd+o(1)}.
\ee

Next we introduce a parameter $N$ which controls the number of terms in $\cA$ that are non-zero,
setting
%. Precisely, we %require that
%define
\be\label{eq:NDef}
N%\cX^{\tau\mu}
\ := \ 
\max\{n\ge1\ : \ a_n\neq0\}
.
\ee
%In practice, this parameter compares to $T$ as follows. 
Since $\by$ is treated as fixed and $\g\in\G$ is of size $T$, we have roughly that $|f(\by\cdot\g)|\le N$, where
\be\label{eq:Nis}
N\ =\ T^{\deg(f)+o(1)}.
\ee

For a square-free integer $q\ge1$ called the {\it level}, we will need to understand the distribution of the sequence $\cA$ along multiples of $q$. To this end, we define
\be\label{eq:cAqDef}
|\cA_{q}|\ := \ \sum_{n\equiv 0(q)}a_{n}(T).
\ee
The following key % ubiquitous
theorem
is used to determine for how large we can take the level and still prove equi-distribution.
\begin{thm}\label{thm:level1}
Let $(\cO,f)$ be as in Examples $A$--$D$, with $\G$ having critical exponent $\gd>1/2$ and spectral gap $\gt<\gd$. For any square-free integer $q\ge1$, we have the estimate % approximation % of the form
\be\label{eq:cAqIs}
|\cA_{q}|
\
=
\
\gw(q)\cdot|\cA|
\ +\ 
O
\left(q\cdot T^{\gt}\,(qT)^{o(1)}\right)
,
\ee
%with the following properties. 
%Assume that 
where
$\gw(q)$ % is
is a ``local density''  function with the following properties.
It is a multiplicative function on  square-free $q$'s
with
\begin{enumerate}
\item
 $\gw(1)=1$,
\item
 for all primes $p\ge2$,
 \be\label{eq:gwp1}
 0\ \le\  \gw(p)\ <\ 1,
 \ee
 and 
 \item
there are constants $K\ge2$ and $\gk\ge1$ so that
we have the local density bound
\be\label{eq:locDens}
\prod_{{z_1\le p \le z}}
\frac1{1-\gw(p)}
\ \le  \
 \left(\frac{\log z}{\log z_1}\right)^\gk
\left(1+{K\over\log z_1}\right)
\ee
for any $2\le z_1<z$.
 \end{enumerate}
\end{thm}

\begin{rmk}
Versions of Theorem \ref{thm:level1} are proved in 
 {\cite[Proposition 4.3]{Kontorovich2009}},  \cite[\S5.2]{KontorovichOh2012}, and  \cite[Theorem 2.1]{LiuSarnak2010} for 
Examples $A$--$D$, but we repeat a sketch of the proof below, as it will be relevant to us later.
\end{rmk}

\begin{rmk}\label{rmk:dim}
%the %rough
 %content
One can 
 interpret 
 \eqref{eq:locDens} as insisting that
the local density at primes be roughly
 \be\label{eq:gwqBig}
\gw(p)\ \approx\ {\gk\over p},
\ee
at least on average, see Lemma \ref{lem:gwp}.
The number $\gk$ appearing in \eqref{eq:locDens} is called the ``sieve dimension'' for $\cA$; see \eqref{eq:gkPyth}. 
Note that $\gk$ is not uniquely defined by \eqref{eq:locDens}, as any larger value also satisfies \eqref{eq:locDens}; in practice one typically takes the least allowable value. 
\end{rmk}

\pf[Sketch of Proof]
To prove Theorem \ref{thm:level1}, 
we first insert the definition \eqref{eq:anTIs}
into \eqref{eq:cAqDef}:
$$
|\cA_{q}| \  =\  
\sum_{n\equiv0(q)}a_{n}(T)
\ = \
\sum_{\g\in\G}
\gU_{T}(\g)\cdot\bo_{\{f(\by\cdot\g)\equiv0(\mod q)\}}
.
$$

The first most basic Affine Sieve observation is that the condition
 \be\label{eq:fgq}
f(\by\cdot \g)\equiv0(\mod q)
\ee 
can be captured by breaking the sum according to the residue of $\g$ mod $q$. In other words, we can decompose 
\be\label{eq:Gdecomp1}
\G\ \cong\ \G(q)\times(\G(q)\bk \G)
.
\ee
%%%%%%%%%
\begin{comment}
\end{comment}
%%%%%%%%%%%%
%
Using this decomposition and following the procedure below, one would obtain \eqref{eq:cAqIs} with the worse error term $O(q^{2}T^{\gt})$, ignoring $o(1)$'s. This would lead to \eqref{eq:ga1} being replaced by the exponent of distribution $\ga=(\gd-\gt)/(3\deg(f))$.

Instead,  what is done in \cite{Kontorovich2009, LiuSarnak2010, KontorovichOh2012} is to capture the condition 
\eqref{eq:fgq}
by decomposing $\by\cdot\g$ (rather than just $\g$) into residue classes $\mod q$. To this end, let $\G_{\by}(q)$ be the stabilizer group of $\by(\mod q)$, that is, define
$$
\G_{\by}(q)\ := \ 
\{\g\in\G \ : \ \by\cdot\g\equiv\by (\mod q)\},
$$
and write $\g\in\G$ uniquely as 
$$
\g\ =\ \g_{1}\g_{0},
$$
with $\g_{1}\in\G_{\by}(q)$ and $\g_{0}\in\G_{\by}(q)\bk\G$. 
Then since $\by\g_{1}\equiv \by(\mod q)$, we have that
\be\label{eq:fbyg1}
f(\by\cdot\g) 
\ =\ 
f(\by\cdot\g_{1}\g_{0}) 
\ \equiv\ 
f(\by\cdot\g_{0}) 
\quad (\mod q)
.
\ee

Hence applying \eqref{eq:gUTAP} with $\Xi(q)=\G_{\by}(q)$, we have 
\bea\nonumber
|\cA_{q}| 
&=& 
\sum_{\g_{0}\in\G_{\by}(q)\bk\G}
\sum_{\g_{1}\in\G_{\by}(q)}
\gU_{T}(\g_{1}\g_{0})
\cdot\bo_{\{f(\by\cdot\g_{1}\g_{0}) \equiv0(\mod q)\}}
\\ \nonumber
&=&
\sum_{\g_{0}\in\G_{\by}(q)\bk\G}
\bo_{\{f(\by\cdot\g_{0}) \equiv0(\mod q)\}}
\left[
\sum_{\g_{1}\in\G_{\by}(q)}
\gU_{T}(\g_{1}\g_{0})
\right]
\\ \label{eq:cAqIs0}
&\overset{\eqref{eq:gUTAP}}{=}&
\sum_{\g_{0}\in\G_{\by}(q)\bk\G\atop
f(\by\cdot\g_{0}) \equiv0(\mod q)}
\left[
\frac{1}{[\G:\G_{\by}(q)]}
|\cA|
+O(T^{\gt+o(1)})
\right]
\\ \label{eq:cAqIs1}
&=&
\frac{\cC(\G_{\by}(q);f)}{[\G:\G_{\by}(q)]}
|\cA|
\ + \ O\left(\cC(\G_{\by}(q);f)\cdot T^{\gt+o(1)}\right)
.
\eea
Here we have defined
\be\label{eq:cCdef}
\cC(\Xi(q);f) \ := \
\#\{\g_{0}\in\Xi(q)\bk\G \ : \
f(\by\cdot\g_{0}) \equiv0(\mod q)\}
,
\ee
where $\Xi(q)$ is  any group with $\G(q)\leqq\Xi(q)\leqq\G$, for which the
 %counting statement 
above makes sense, that is, whenever the condition $f(\by\cdot\g_{0})\equiv0(\mod q)$ is left-$\Xi(q)$ invariant. 
%\begin{rmk}\label{rmk:wellDef}

%\end{rmk}
%
Now we can set the local density function to be
\be\label{eq:gwqIs}
\gw(q)\ :=\ \frac{\cC(\G_{\by}(q);f)}{[\G:\G_{\by}(q)]},
\ee 
whence we have a decomposition of the form \eqref{eq:cAqIs}. 

It is straightforward to compute that the index 
\be\label{eq:index1}
[\G:\G_{\by}(q)]=q^{2+o(1)},
\ee 
and moreover that, very roughly,
\be\label{eq:cCBnd1}
\cC(\G_{\by}(q);f)\ < \ q^{1+o(1)}
.
\ee
Inserting \eqref{eq:cCBnd1} into the error term of \eqref{eq:cAqIs1} confirms the error term in \eqref{eq:cAqIs}.
%, and hence the exponent of distribution \eqref{eq:ga1}. 

It remains to verify the properties of $\gw$. The condition (1), that is, that $\gw(1)=1$, is clear, and multiplicativity 
%is not hard to establish (
%in the thin orbit cases, these
 follows from Strong Approximation and Goursat's Lemma. %, see \cite{Kontorovich2009, KontorovichOh2012}). 
It follows from the strong primitivity assumption that $\gw(p)<1$ for all primes. 
Verification of the key property \eqref{eq:locDens} is postponed to the next Lemma, whence the proof of Theorem \ref{thm:level1} is complete.
\epf

The following Lemma verifies \eqref{eq:gwqBig}, from which the local density bound \eqref{eq:locDens} follows by classical methods.

\begin{lem}\label{lem:gwp}
For primes $p$ sufficiently large, we have the following estimates on $\gw(p)$. 
In the ``thin Pythagorean'' cases, we have that (see  \cite[Lemma 5.4]{KontorovichOh2012}):
$$
\gw(p)=\twocase{}
{{2\over p+1},}{if $p\equiv1(\mod 4)$,}
{0,}{if $p\equiv3(\mod 4)$,}
\qquad\qquad\text{
%if $f=f_{\sH}$,
for Example $A$,
}
$$
$$
\gw(p)=
{{4\over p+1}}
,
\qquad\qquad\text{
%if $f=f_{\sA}$,
for Example $B$,}
$$
and\footnote{This corrects a typo in \cite[(5.6)]{KontorovichOh2012}.}
$$
\gw(p)=\twocase{}
{{6\over p+1},}{if $p\equiv1(\mod 4)$,}
{{4\over p+1},}{if $p\equiv3(\mod 4)$,}
\qquad\qquad\text{
%if $f=f_{\sC}$.
for Example $C$.}
$$

In the ``anisotropic'' case, we have \cite[(6.4)]{LiuSarnak2010} that
$$
\gw(p)=
\frac3p+O\left(\frac1{p^{2}}\right)
,\qquad\qquad\text{
%if $f=f_{\sC}$.
for Example $D$.}
$$
In particular, \eqref{eq:locDens} holds with
\be\label{eq:gkIs}
\gk=
\fourcase
{1,}{in Example $A$,}
{4,}{in Example $B$,}
{5,}{in Example $C$,}
{3,}{in Example $D$.}
\ee
\end{lem}
%In converting the Lemma to 
% the local density bound \eqref{eq:locDens},
% we find that
%Examples $A$--$D$ have sieve dimensions $\gk=1,4,5$ and $3$, respectively, as claimed.
\begin{rmk}\label{rmk:IsoAni}
It is only here in the local density estimate that the sieve can distinguish the sieve dimensions $\gk$ for the product of coordinates in the isotropic Example $C$ and the anisotropic Example $D$ (see  Remark \ref{rmk:IsoAni0}). 
The form $\cQ$ being isotropic is equivalent to the cone $\cQ=0$ (and other level sets) being parametrizable by non-constant polynomial maps. In particular, if $\cQ$ is isotropic, then there exist rational binary quadratic forms $G_{1},G_{2},G_{3}$ so that
$$
\cQ(G_{1}(c,d),G_{2}(c,d),G_{3}(c,d))=0.
$$
If the $G_{j}$ are reducible, then the product of coordinates $f(\bx)=xyz$  can be the product of more than $3$ irreducible factors, and this is exactly what happens in the Pythagorean case. 
On the other hand, no such  parametrization exists if $\cQ$ is anisotropic, whence the product of three coordinates always has sieve dimension $\gk=3$.

 To make this completely concrete for the form $\cQ=x^{2}+y^{2}-z^{2}$,
recall the ancient parametrization of Pythagorean triples $\bx=(x,y,z)$ with $y$ even as
$$%\beann
%x&=&c^{2}-d^{2},\\
%y&=&2cd,\\
%z&=&c^{2}+d^{2}.
\threecase
{x\ =\ G_{1}(c,d)\ =\ c^{2}-d^{2},}{}
{y\ =\ G_{2}(c,d)\ =\ 2cd,}{}
{z\ = \ G_{3}(c,d)\ =\ c^{2}+d^{2}.}{}
$$%\eeann
Both $G_{1}$ and $G_{2}$ factor into products of two linear forms, and so in Example $C$,
$$
f_{\sC}(\bx) \ = \ \frac1{60}xyz \ = \ \frac1{30} (c+d)(c-d)cd(c^{2}+d^{2})
$$
is a product of $\gk=5$ irreducible factors. 
 
 On the other hand, the form 
 $$
 \cQ(\bx)=x^{2}+y^{2}-2z^{2}
 $$ 
 is also isotropic over $\Q$, but the cone $\cQ=0$ has a parametrization
$$%\beann
%x&=&c^{2}-d^{2},\\
%y&=&2cd,\\
%z&=&c^{2}+d^{2}.
\threecase
{x\ =\ G_{1}(c,d)\ =\ c^{2}+2cd-d^{2},}{}
{y\ =\ G_{2}(c,d)\ =\ c^{2}-2cd-d^{2},}{}
{z\ = \ G_{3}(c,d)\ =\ c^{2}+d^{2},}{}
$$%\eeann
in which all three forms $G_{j}$ are irreducible. In this example, the product of coordinates would have  sieve dimension  $\gk=3$. 
%See \cite{Marasingha2010}, where this example 
\end{rmk}

%\begin{rmk}
In light of \eqref{eq:gwqBig} and \eqref{eq:cAsize}, the ``main''  term in the approximation \eqref{eq:cAqIs} is roughly of size  $T^{\gd}/q$, while the ``error'' is % roughly of size 
about
$qT^{\gt}$. Balancing these terms,
we %should 
can
take $q$ 
%in the range
%is a bit less than
almost as large as
$T^{(\gd-\gt)/2}$.
Converting to the parameter $N$ in \eqref{eq:Nis},  we see that the approximation \eqref{eq:cAqIs} is a 
%this gives a 
true asymptotic %only 
whenever
\be\label{eq:lev1}
q\ <\ N^{(\gd-\gt)/(2\deg(f))-\vep},
\ee
for any fixed $\vep>0$.  
%We make this precise below.
%\end{rmk}
For later reference, we record the following estimate, which follows immediately from \eqref{eq:cAqIs}. %, we have the following immediate
\begin{cor}\label{cor}
For any fixed $\vep>0$, there is an $\eta=\eta(\vep)>0,$
 %$\eta\to0$ as $\vep\to0$, 
 so that
\be\label{eq:level}
\sum_{q<N^{\ga-\vep}\atop q \text{ square-free}}
\bigg|
|\cA_{q}|
\
-
\
\gw(q)\cdot|\cA|
\bigg|
\ \ll_{\vep}\
|\cA|^{1-\eta}, %+o(1)},
\ee
as $T\to\infty$,
where
\be\label{eq:ga1}
\ga\ := \ {\gd-\gt\over2\deg(f)}
\ee
is the exponent in \eqref{eq:lev1}.
\end{cor}

 \begin{rmk}
The quantity
 %on the right side of \eqref{eq:lev1} 
 $N^{\ga}$
 is called a {\it level of distribution} for $\cA$, and the exponent
$\ga$ in \eqref{eq:ga1}
is called an {\it exponent of distribution}. 
%The %level or 
%exponent of distribution
This
 is not a quantity intrinsic to $\cA$ but is rather a function of what one can prove about $\cA$. In particular,
any smaller value of $\ga$ is also an exponent  of distribution, but in applications, one wishes to take $\ga$ as large as possible.
\end{rmk}

\begin{rmk}\label{rmk:corrKO}
%\footnote{\label{foot}%
We are correcting here a %^{}minor mistake 
typo
in \cite[(2.23)]{KontorovichOh2012}, where
 %the first lower bound on $\mu$ should be $2\deg(f)/(\gd-\gt)$ instead of $2/(\gd-\gt)$; 
 $\deg(f)$ was omitted from $\ga$ (our $\ga$ is $1/\mu$ in the notation of \cite{KontorovichOh2012}); 
hence the values of $R$ computed in \cite{KontorovichOh2012} are only accurate in the case $f=f_{\sH}$ of Example $A$, see Remark \ref{rmk:KO}.
\end{rmk}

%Note that in other applications, one can take advantage 

\begin{rmk}\label{rmk:avg}
In sieve applications, one only needs the average estimate \eqref{eq:level} and not the estimate for  individual $q$ discussed before \eqref{eq:lev1}.
% such an estimate to hold on average over $q$;
In Example $A$, it is exactly this averaging which is exploited in \cite{BourgainKontorovich2013} to prove  \eqref{eq:BK}. In Examples $B$--$D$, we do not currently know how to exploit this average, and so the level of distribution just follows from the individual estimate \eqref{eq:cAqIs}.
See \cite{Marasingha2010} for some  sharp levels of distribution for non-thin isotropic (and hence parametrizable; cf. Remark \ref{rmk:IsoAni}) orbits, also obtained by exploiting the average on $q$.
\end{rmk}

We now have all the properties we need from the sequence $\cA$. In the next subsection, we recall the high-dimensional weighted sieve used in applications.

\subsection{Diamond-Halberstam-Richert Sieve}\label{sec:sieve}\

Recall that $\cP_{R}$ is the set of $R$-almost primes. Sieve theory  produces an estimate for
$$
\sum_{n\in\cP_{R}}a_{n}(T),
$$
given 
knowledge of the distribution of $\cA$ along arithmetic progressions.
Adapted to our setting, we have the following

\begin{figure}
${\includegraphics[width=2in]{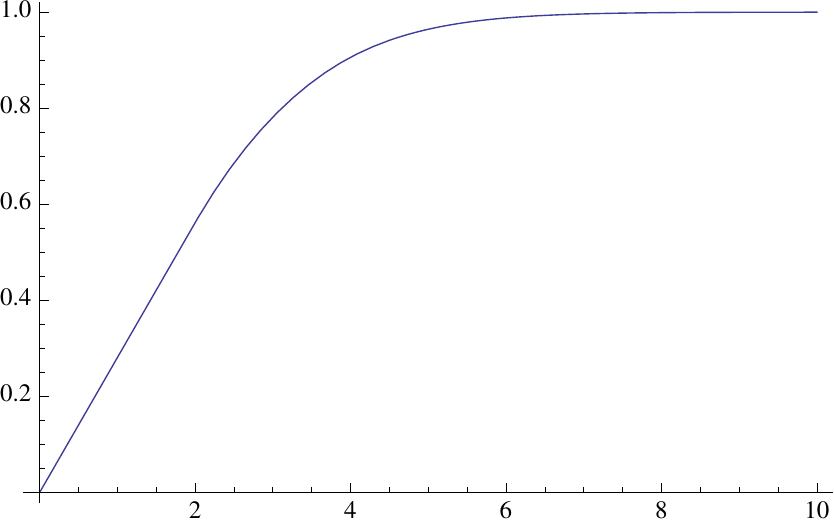}\atop \gs(u)}
{\includegraphics[width=2in]{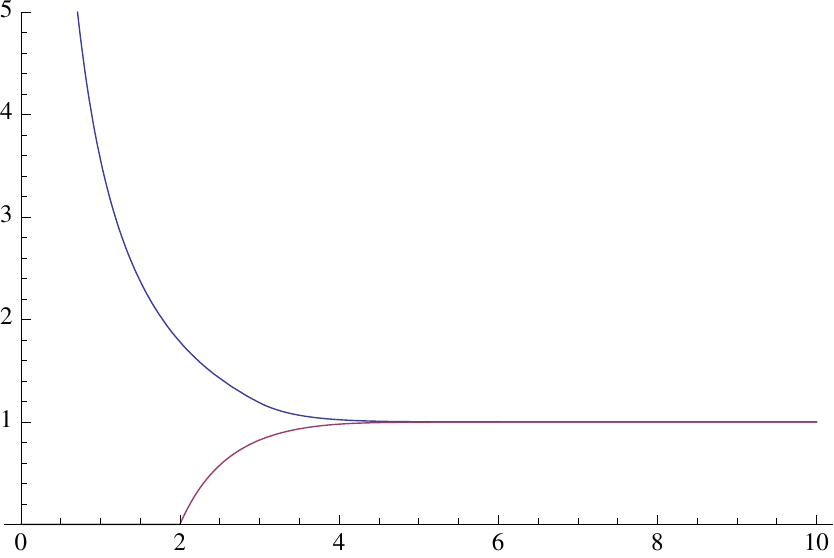}\atop f(u),\ F(u)}$
\caption{Plots  of $\gs(u)$, $f(u)$ and $F(u)$ for $\gk=1$.}
\label{fig:gsfF}
\end{figure}
\begin{Thm}[\cite{DiamondHalberstamRichert1988,DiamondHalberstam1997}]\label{thm:sieve}
Let $\cA$, $N$,
%$z$,
$\gw$, $\gk$, and
$\ga$
% and $\tau
 be as described in
\eqref{eq:anTIs},
 \eqref{eq:NDef},
 \eqref{eq:gwqIs}, 
 \eqref{eq:gkIs},
%\eqref{eq:cAqIs}, 
and \eqref{eq:ga1};
in particular, they satisfy the key conditions \eqref{eq:locDens} and \eqref{eq:level}.
It is convenient to define another parameter
\be\label{eq:tauDef}
\tau \ := \
{\ga \log N\over \log |\cA|}
\ = \ {\ga \cdot \deg(f)\over \gd} +o (1).
\ee

%\begin{enumerate}
%\item
$(i)$
Let $\gs(u)=\gs_{\gk}(u)$ be the continuous solution of the differential-difference problem:
\be
\begin{cases}
{u^{-\gk}\gs(u)=A_{\gk}^{-1},} & \text{for $0<u\le2$, $A_{\gk}=(2e^{\g})^{\gk}\G(\gk+1)$}\\
{(u^{-\gk}\gs(u))'=-\gk u^{-\gk-1} \gs(u-2),}&\text{for $u>2$,}
\end{cases}\ee
where $\g$ is the Euler constant and $\G$ is the Gamma function\footnote{There should be no confusion here with the discrete group $\G$.}. Then
there exist numbers
\be\label{eq:gagb}
\ga_{\gk}\ge\gb_{\gk}\ge2
\ee
so that
 the following simultaneous differential-difference system has continuous solutions $F(u)=F_{\gk}(u)$ and $f(u)=f_{\gk}(u)$ which satisfy
$$
F(u)=1+O(e^{-u}), \quad f(u)=1+O(e^{-u}),
$$
 and $F$ (resp. $f$) decreases (resp. increases) monotonically towards $1$ as $u\to\infty$:
\be
\begin{cases}
{F(u)=1/\gs(u),} & \text{for $0<u\le\ga_{\gk},$}\\
{f(u)=0,} & \text{for $0<u\le\gb_{\gk},$}\\
{(u F(u))'= f(u-1),} &\text{for $u>\ga_{\gk},$}\\
{(u f(u))'= F(u-1),} &\text {for $u>\gb_{\gk}.$}
\end{cases}\ee
See Figure \ref{fig:gsfF} for plots of $\gs,$ $f$ and $F$ in the case $\gk=1$. 

%\item
$(ii)$
For any two real numbers $u$ and $v$ with
\be\label{eq:taugb}
 \tau^{-1}<u\le v,\qquad \gb_{\gk}<\tau v,
\ee
and assuming that
\be\label{Ris}
R\ > \ {\tau u\over\ga} - 1 +
 {\gk \over f(\tau v)} \int_1^{v/u}F(\tau v-s)\left(1-\frac uv s\right){ds\over s}
 ,
\ee
we have
\be\label{eq:anTLower}
\sum_{n\in\cP_{R}}a_{n}(T)\  \gg\  |\cA| \prod_{p<N^{\ga}}(1-\gw(p))
\ \gg \ 
{|\cA|\over (\log T)^{\gk}}
.
\ee
%\end{enumerate}
\end{Thm}

\begin{table}
\begin{tabular}{|r|c|c|c|c|}
\hline
$\gk$&1&3&4&5\\ \hline

%$\ga_{\gk}$ & &&&\\ \hline

$\gb_{\gk}$ & \ \ 2 \ \  &6.6408\dots&9.0722\dots&11.5347\dots
\\
\hline
\end{tabular}
\vskip.1in
 \caption{Values of %$\ga_{\gk}$ and 
 $\gb_{\gk}$ for $\gk=1,3,4,5$.}
\label{table1}
\end{table}

\begin{rmk}\label{rmk:Zariski}
The statements in % Theorems \ref{thm:KO}, \ref{thm:Pyth}, \ref{thm:LS}, and \ref{thm:Aniso} 
Examples $A$--$D$
on the  Zariski density of $\cO(f,R)$ are then proved  easily from the archimedean bounds in \eqref{eq:anTLower}; see, e.g., the proof of \cite[Corollary 2.3]{LiuSarnak2010}.
\end{rmk}

\begin{rmk}
The sieve dimensions relevant to us are $\gk=1,3,4,$ and $5$, and we will need the corresponding values of the  constant $\gb_{\gk}$ for \eqref{eq:taugb}.
%(and $\ga_{\gk}$, though this variable is less important) have
 %been
 These
 are computed in \cite[p. 345]{DiamondHalberstam1997}, and we reproduce them in Table \ref{table1}. 
\end{rmk}

While the expression on the right hand side of \eqref{Ris} is unwieldy, it can certainly be estimated by  one's favorite software package.
That said,  the following simplification is quite effective in practice (see \cite[(6.15)]{LiuSarnak2010}):  for any $0<\gz<\gb_{\gk}$, the expression  is maximized by any value of
\be\label{eq:mIs}
m_{\ga,\gk}(\gz)\ := \
{1\over\ga}\left(1+\gz-{\gz\over\gb_{\gk}}\right) - 1 +
(\gk+\gz)\log{\gb_{\gk}\over \gz}-\gk+\gz{\gk\over\gb_{\gk}}
 .
\ee
%where we have set
%$$
%\tau u=1+\gz-{\gz\over\gb_{\gk}},\qquad
%\tau v={\gb_{\gk}\over\gz}+\gb_{\gk}-1.
%$$

\subsection{Proofs of Theorems \ref{thm:KO} and \ref{thm:LS}}\label{sec:affStd}\

It remains to insert the specific values of $\ga$, $\gk$, and $\tau$, and compute the resulting values of $R$ for each of our examples. 

\subsubsection{Example $A$}

To obtain as small a value of $R$ as possible, we take $\gd$ as large as possible, that is, near $1$, to take advantage of Gamburd's $\gt=5/6$ gap in Theorem \ref{thm:spec}. At first, we just set $\gd=1$. The degree of $f=f_{\sH}$ is $\deg(f)=1$, so the exponent of distribution \eqref{eq:ga1} is
$$
\ga \ = \ {1-5/6\over2}\ =\ \frac1{12},
$$
and the sieve dimension is $\gk=1$. With these values of $\ga$ and $\gk$, the minimal value of $m(\gz)$ in \eqref{eq:mIs} is $m(0.12)=13.93$, leading to the bound $R_{0}(\cO,f_{\sH})\le 14$ for the saturation number. Letting $\gd$ be slightly less than $1$, we can still ensure that $\ga$ is large enough that the minimal value of $m(\gz)$ is $<14$. 

But in fact,  better methods are known to estimate $R$-values for linear (that is, dimension $\gk=1$) sieve problems using essentially identical assumptions, see, e.g., Richert's weights in \cite[\S25.3]{FriedlanderIwaniecBook}. From the exponent of distribution $\ga_{A}=1/12$, these produce $R=13$-almost primes, with room to 
%decrease 
perturb
$\gd$ to a little below $1$. 
This is the $R$ value we stated in Theorem \ref{thm:KO} for  Example $A$. Regardless, none of these values are
 %interesting
 relevant
  anymore, in light of \eqref{eq:BK}.

\subsubsection{Example $B$}

Because the Pythagorean form $\cQ(\bx)=x^{2}+y^{2}-z^{2}$ 
is isotropic with reducible parametrizing forms (see Remark \ref{rmk:IsoAni}), the sieve dimension for the ``area'' function $f(\bx)=f_{\sA}(\bx)=\frac1{12}xy$ is $\gk=4$ (rather than $\gk=2$). The degree is $\deg(f)=2$. As above, we begin by taking $\gd=1$ and using Gamburd's gap $\gt=5/6$. This gives the exponent of distribution
$$
\ga  \ = \ {1-5/6\over2\cdot 2}\ =\ \frac1{24}.
$$
Optimizing $m(\gz)$ with these values gives $m(0.16)=39.28$. Again, letting $\gd$ be slightly below $1$ still recovers the value $R=40$, as claimed in Theorem \ref{thm:KO}. 

\subsubsection{Example $C$}

For $f=f_{\sC}$, the degree is $\deg(f)=3$ and sieve dimension is $\gk=5$. Again we take $\gd=1$ and $\gt=5/6$, giving the exponent of distribution
$$
\ga  \ = \ {1-5/6\over2\cdot 3}\ =\ \frac1{36}.
$$
Now optimizing $m(\gz)$ gives $m(0.136)=57.3$. For $\gd$ slightly below $1$, we  still recover $R=58$.

\subsubsection{Example $D$}

In this non-thin anisotropic example, we have $\gd=1$, $\deg(f)=3$, and sieve dimension $\gk=3$. Using the Kim-Sarnak spectral gap $\gt=39/64$ in Theorem \ref{thm:spec}, we obtain the exponent of distribution
$$
\ga \ = \ {1-39/64 \over 2\cdot 3}\ \approx\  {1\over15.36}.
$$
Optimizing $m(\gz)$ gives $m(0.186)=25.26$, giving the claimed value $R=26$. Assuming the Selberg Eigenvalue Conjecture, we can take $\gt=1/2$ and 
$$
\ga \ = \ {1-1/2 \over 2\cdot 3} \ =\ {1\over12}.
$$
Then $m(\gz)$ is optimized at $m(0.23)=21.3$, giving $R=22$, as claimed in Theorem \ref{thm:LS}.
\\

These are the  values of $R$ produced  in \cite{Kontorovich2009}, \cite{LiuSarnak2010}, and  \cite{KontorovichOh2012}. In the next section, we make one further simple observation, which has the effect of doubling
 %the value of
  the exponent of distribution over that in \eqref{eq:ga1}. 
Using the same methods as here, we then conclude Theorems \ref{thm:Pyth} and \ref{thm:Aniso}.

\newpage

\section{Proofs of  Theorems \ref{thm:Pyth} and \ref{thm:Aniso}}\label{sec:pfs}

We keep all the same notation from the previous section, first describing %the 
our
initial
aim 
%of this project
 in rough terms, before explaining our new observation. 

\subsection{Initial Idea}\

The 
goal 
of this project was to try to 
improve the level of distribution 
by exploiting the $\g_{0}$ sum
  in
 \eqref{eq:cAqIs0}, which was estimated trivially to arrive at \eqref{eq:cAqIs1}. Of course this requires us to keep track of all the lower order terms in \eqref{eq:gUTAP}, rather then estimating them in absolute value. We proceed as follows.

We will want $\Xi(q)\bk \G$ to be a group (i.e. $\Xi(q)$ to be normal in $\G$), so return to the decomposition \eqref{eq:Gdecomp1}; that is, we set $\Xi(q)=\G(q)$, rather than $\Xi(q)=\G_{\by}(q)$. 
(So the length of the $\g_{0}$ sum in  \eqref{eq:cAqIs0} is now about $q^{2}$ instead of $q$, but we hope to recover this loss and more.)
Assume for simplicity that $\G(q)\bk\G\cong\PSL_{2}(q)$ and that $q$ is prime.
%In fact,
% 
% \begin{rmk}\label{rmk:Project}
%Then 
The space $L^{2}(\G(q)\bk G)$ carries not only a right (regular) $G$-action, but also a left (Hecke-like) $\G(q)\bk\G$-action. %, which of course commute. 
Decomposing %first 
with respect to the latter action,
the estimate \eqref{eq:gUTAP} 
%is the result of 
can be % interpreted 
%seen
%as
%is 
obtained from
an
 %decomposition 
expansion
 of the form  (see the discussion after \cite[(2.12)]{BourgainKontorovich2013})
\be\label{eq:gUTAP2}
\sum_{\g_{1}\in\G(q)}\gU_{T}(\g_{1}\g_{0}) \ = \
\sum_{\rho\in \widehat{
%\G(q)\bk\G
\PSL_{2}(q)
}}
\cM_{\rho}(T,q;\g_{0})
,
\ee
where $\rho$ ranges over irreducible unitary representations of the finite group $\
%Xi(q)\bk\G
\PSL_{2}(q)
$, and $\cM_{\rho}$ is the contribution coming from  $\rho$. The first term in \eqref{eq:gUTAP} comes from $\rho=\bo$, that is, the trivial representation;  
the other $\rho$'s come from new spectrum, and \eqref{eq:gUTAP}
is obtained by
%and thus are
 controlling these terms in totality by the spectral gap. %, resulting in
%. %The relevance of this remark is made clear in \S\ref{sec:pfs}.
%\end{rmk}

 Instead of estimating %directly 
 these
 the error terms in absolute value, %as done in \eqref{eq:gUTAP}, 
 we 
 will 
 want
 %ed 
 to capitalize on the full decomposition \eqref{eq:gUTAP2}. So we insert it into the analogue of \eqref{eq:cAqIs0}, capture the condition $f(\by\cdot \g_{0})\equiv0(q)$ by abelian harmonic analysis, and carry out the $\g_{0}\in\PSL_{2}(q)$ sum on each irreducible. Expanding out the terms, one faces the following problem, which seems to be new (see the somewhat related questions arising in \cite{SotoAndrade1987, Katz1993}): Given an irreducible unitary representation $(\rho,V)$ of a finite non-abelian group, e.g. $\PSL_{2}(\bbF_{q})$,
 %, say of $\SL_{2}(\bbF_{p})$, 
 an additive character $\psi$ on $\bbF_{q}$, and a polynomial $f$ on $\mattwos abcd\in\PSL_{2}$, capture cancellation in the matrix coefficients of $\rho$ against the character of the polynomial; that is, give a non-trivial estimate for a sum of the form
 $$
%V^{2}\ni (v,w)\ \mapsto\ 
\sum_{\g\in\PSL_{2}(\bbF_{q})}
\<\rho(\g).v,w\>
\cdot 
\psi(f(\g))
,
 $$
for $v,w\in V$.
While the initial aim of the project was somewhat sophisticated, 
after computing several explicit examples of the above type, 
we stumbled upon a  completely elementary observation that had been previously overlooked.
Its effect, in our applications, is to make the $\g_{0}$ sum in the analogue of \eqref{eq:cAqIs0} have length $q^{\vep}$, rather than $q$, leading to a level of distribution twice as large as before. 
So while the above general problem is still interesting and may have other applications,  in the end it is of no consequence to our current results.

%\newpage

\subsection{The observation}\label{sec:obs}
\

The key new observation is that one can use a larger group than $\G_{\by}(q)$ in capturing the condition \eqref{eq:fgq}. To this end, we introduce the group
$\G_{\<\by\>}(q) $
 which stabilizes the linear span of $\by$ mod $q$. That is, we define
\beann
\G_{\<\by\>}(q) & := &
\{
\g\in\G
\ : \
\by\cdot \g\in\<\by\>(\mod q)
\}
\\
&=&
\{
\g\in\G
\ : \
\exists a\in(\Z/q\Z)^{\times}\text{ with }\by\cdot \g\equiv a\by(\mod q)
\}
.
\eeann
Clearly $$\G(q)\leqq\G_{\<\by\>}(q)\leqq\G.$$
%The key o
%Observe
%aetion is 
Note
that, because the functions $f$ in all the Examples $A$--$D$ are homogeneous, we have
$$
f(\by\cdot\g_{1}\g_{0}) 
\ \equiv\ 
a^{\deg(f)}f(\by\cdot\g_{0}) 
\quad (\mod q)
,\qquad
\text{ for some }a\in(\Z/q\Z)^{\times},
$$
whenever $\g_{1}\in\G_{\<\by\>}(q)$.  Hence we can replace \eqref{eq:fbyg1} by the fact that 
\be\label{eq:obs}
f(\by\cdot\g_{1}\g_{0})\equiv0(\mod q)\quad\text{ if and only if }\quad f(\by\cdot\g_{0})\equiv0(\mod q).
\ee 

\begin{rmk}
We emphasize here that it is not the cone $\cQ=0$  which takes advantage of this homogeneity (since we also consider other level sets, see \eqref{eq:levelSet}), but
rather
 the   sieve,
which only asks for the distribution of $a_{n}(T)$ on {\it multiples} of $q$, see \eqref{eq:cAqDef}. If for other applications one wants to capture residue classes other than $0(\mod q)$, then the homogeneity of $f$ will not help.
\end{rmk}

Beyond this simple observation, we proceed exactly as described in \S\ref{sec:2}.

\subsection{The Proofs}\

%We replace the key Theorem \ref{thm:level1} with
\begin{thm}\label{thm:level3}
Theorem \ref{thm:level1} holds exactly as stated, but with \eqref{eq:cAqIs} replaced by
\be\label{eq:cAqIsN}
|\cA_{q}|
\
=
\
\gw(q)\cdot|\cA|
\ +\ 
O
\left( T^{\gt}\,(qT)^{o(1)}\right)
.
\ee
\end{thm}

\begin{rmk}\label{rmk:new}
The only difference to notice is that the error term in \eqref{eq:cAqIs3} is essentially $T^{\gt}$, rather than $qT^{\gt}$ in \eqref{eq:cAqIs}; that is, we have recovered a power of $q$ which was lost in the previous approach. This explains our comment in \S\ref{sec:new}.
\end{rmk}

\pf[Sketch of Proof]
We start with the same definition of $a_{n}(T)$ as in \eqref{eq:anTIs}.
Replacing the decomposition \eqref{eq:Gdecomp1} with
\be\label{eq:Gdecomp3}
\G\ \cong\ \G_{\<\by\>}(q)\times(\G_{\<\by\>}(q)\bk \G),
\ee
 we now write
\bea\nonumber
|\cA_{q}| 
&=& 
\sum_{\g_{0}\in\G_{\<\by\>}(q)\bk\G}
\sum_{\g_{1}\in\G_{\<\by\>}(q)}
\gU_{T}(\g_{1}\g_{0})\cdot
\bo_{\{f(\by\cdot\g_{1}\g_{0}) \equiv0(\mod q)\}}
\\ \nonumber
&\overset{\eqref{eq:obs}}{=}&
\sum_{\g_{0}\in\G_{\<\by\>}(q)\bk\G}
\bo_{\{f(\by\cdot\g_{0}) \equiv0(\mod q)\}}
\left[
\sum_{\g_{1}\in\G_{\<\by\>}(q)}
\gU_{T}(\g_{1}\g_{0})
\right]
\\ \nonumber
&\overset{\eqref{eq:gUTAP}}{=}&
\sum_{\g_{0}\in\G_{\<\by\>}(q)\bk\G\atop
f(\by\cdot\g_{0}) \equiv0(\mod q)}
\left[
\frac{1}{[\G:\G_{\<\by\>}(q)]}
|\cA|
+O(T^{\gt})
\right]
\\ \label{eq:cAqIs3}
&=&
\frac{\cC(\G_{\<\by\>}(q);f)}{[\G:\G_{\<\by\>}(q)]}
|\cA|
\ + \ O\left(\cC(\G_{\<\by\>}(q);f)\cdot T^{\gt}\right)
,
\eea
where we applied \eqref{eq:gUTAP} with $\Xi(q)=\G_{\<\by\>}(q)$, and used the definition \eqref{eq:cCdef}. Note we are allowed to use $\Xi(q)=\G_{\<\by\>}(q)$ in \eqref{eq:cCdef}; indeed, the observation \eqref{eq:obs}  says precisely  that %the condition 
$f(\by\cdot\g_{0})\equiv0(\mod q)$  (as a condition on $\g_{0}$) is left-$\G_{\<\by\>}(q)$ invariant. 

The new ``local density'' function 
\be\label{eq:gwqNew}
\gw(q)\ :=\ \frac{\cC(\G_{\<\by\>}(q);f)}{[\G:\G_{\<\by\>}(q)]}
\ee 
is then actually the same as that in \eqref{eq:gwqIs}. Indeed, just fix $q$ and take $T\to\infty$, comparing \eqref{eq:cAqIs3} with \eqref{eq:cAqIs1}.
Thus the sieve dimensions are the same as before.

On the other hand, the index $[\G:\G_{\<\by\>}(q)]$ is now of size $q^{1+o(1)}$ instead of \eqref{eq:index1}.
%. Using similar  techniques to those leading to Lemma \ref{lem:gwp}, one then establishes  the bound
Thus comparing \eqref{eq:gwqNew} to Lemma \ref{lem:gwp} gives
\be\label{eq:cCbnd3}
\cC(\G_{\<\by\>}(q);f)<q^{o(1)}
\ee
instead of \eqref{eq:cCBnd1}. Inserting \eqref{eq:cCbnd3} into \eqref{eq:cAqIs3} gives \eqref{eq:cAqIsN}, as claimed.
\epf

Then we obtain the same Corollary \ref{cor} but with the exponent of distribution
\be\label{eq:gaIs3}
\ga\  =\  {\gd-\gt\over \deg(f)},
\ee
instead of \eqref{eq:ga1}. That is, the effect of the simple observation \eqref{eq:obs} is to 
 double the exponent of distribution. 

With all the other ingredients in place, it remains to estimate the new values of $R$.

\subsubsection{Example $A$}\label{sec:ExA}

As before, we start by taking $\gd=1$ with Gamburd's $\gt=5/6$ spectral gap. The sieve dimension and degree are both $\gk=\deg(f)=1$.
Inserting these values into \eqref{eq:gaIs3} gives the exponent of distribution 
$$
\ga\  =\  {1-5/6\over 1}\ = \ \frac16.
$$
Linear sieve methods then produce the value $R=7$, with room to allow $\gd$ a little below $1$;
see
Remark \ref{rmk:ExA}.
%Of course, if we allow here  the full group $\G=\SO_{\cQ}(\Z)$, then we obtain all Pythagorean triples. Fermat showed that infinitely many Pythagorean triples have prime hypotenuse

\subsubsection{Example $B$}

Again we take $\gd=1$, $\gt=5/6$, and $\deg(f)=2$. The exponent of distribution is
$$
\ga\  =\  {1-5/6\over 2}\ = \ \frac1{12}
$$
for this dimension $\gk=4$ problem.
Optimizing the function $m(\gz)$ in \eqref{eq:mIs} gives $m(0.295)=24.99,$ or $R=25$. 

\subsubsection{Example $C$}

We set $\gd=1$, with $\gt=5/6$ and $\deg(f)=3$. The exponent of distribution is then
$$
\ga\  =\  {1-5/6\over 3}\ = \ \frac1{18}
$$
for a sieve of dimension $\gk=5$.
Optimizing $m(\gz)$ gives $m(0.25)=36.3,$ or $R=37$.  This completes the proof of Theorem \ref{thm:Pyth}.

\subsubsection{Example $D$}

Unconditionally, we have $\gd=1$, and the Kim-Sarnak gap $\gt=39/64$ with $\deg(f)=3$. The exponent of distribution is then
$$
\ga\  =\  {1-39/64\over 3}\ \approx \ \frac1{7.7}.
$$
The sieve dimension is $\gk=3$, and 
optimizing $m(\gz)$ gives $m(0.33)=15.9,$ or $R=16$.  

Assuming the Selberg Eigenvalue Conjecture, we can take $\gt=1/2$ with exponent
$$
\ga\  =\  {1-1/2\over 3}\ = \ \frac1{6}.
$$
Now
optimizing $m(\gz)$ gives $m(0.4)=13.7,$ or $R=14$.  

This completes the proof of Theorem \ref{thm:Aniso}.

%\begin{comment}
%\end{comment}

\newpage

%%%%%%%%%%%%%%%%%%%%%%%%%%%%%%%%%%%%%%%%

\bibliographystyle{alpha}

\bibliography{AKbibliog}

\end{document}